\documentclass{elsarticle}

\usepackage{amsmath}
\usepackage{bm}
\usepackage{tikz}
\usepackage{pgfplots}

\usepackage{hyperref}

\newproof{pf}{Proof}

\journal{arXiv.org} 

\graphicspath{{pictures/}}

\begin{document}

\begin{frontmatter}

\title{Numerical solution of boundary value problems for the eikonal equation in an anisotropic medium}

\author[nsi]{Alexander~G.~Churbanov}
\ead{achur@ibrae.ac.ru}

\author[nsi,univ]{Petr N. Vabishchevich\corref{cor}}
\ead{vabishchevich@gmail.com}

\address[nsi]{Nuclear Safety Institute, Russian Academy of Sciences, 52, B. Tulskaya, Moscow, Russia}
\address[univ]{North-Eastern Federal University, 58, Belinskogo, Yakutsk, Russia}

\cortext[cor]{Corresponding author}

\begin{abstract}
A Dirichlet problem is considered for the eikonal equation in an anisotropic medium.
The nonlinear boundary value problem (BVP) formulated in the present work is the limit of the 
diffusion--reaction problem with a diffusion parameter tending to zero.
To solve numerically the singularly perturbed diffusion--reaction problem,
monotone approximations are employed. Numerical examples are presented for a two-dimensional
BVP for the eikonal equation in an anisotropic medium. The standard piecewise-linear
finite-element approximation in space is used in computations.
\end{abstract}


\end{frontmatter}

\section{Introduction}

Many applied problems lead to the need of solving a BVP for the eikonal equation.
First of all, this nonlinear partial differential equation is used to simulate wave propagation 
in the approximation of geometric optics \cite{Born,Kravtsov}.
In computational fluid dynamics, image processing and computer graphics (see, for example, 
\cite{sethian1999level,Gilles}),
the solution of BVPs for the eikonal equation is associated with calculating the nearest distance 
to boundaries of a computational domain.

The eikonal equation is a typical example of steady-state Hamilton--Jacobi equations.
The issues of the existence and uniqueness of the solution for boundary value problems 
for such equations are considered, e.g., in \cite{kruvzkov1975generalized,lions1982generalized}.
To solve numerically BVPs for the eikonal equation, the standard approaches are used, 
which are based on using difference methods on rectangular grids or finite-element/finite-volume 
approximations on general irregular grids. In this approach, the main attention is paid to problems of nonlinearity.

A boundary value problem is formulated in the following way.
The function $u(\bm x)$ is defined as the solution of the equation
\begin{equation}\label{1}
 |\nabla u|^2 = 1,
 \quad \bm x \in \Omega 
\end{equation} 
in a domain $\Omega$ with the specified boundary conditions
\begin{equation}\label{2}
 u(\bm x) = 0,
 \quad \bm x \in \partial \Omega . 
\end{equation}

Computational algorithms for solving BVPs for the eikonal equation can be divided into two classes.

Marching methods (the first class of algorithms) are the most widely used.
They are based on the hyperbolic nature of the eikonal equation.
In this case, the desired solution of the problem (\ref{1}), (\ref{2})
is obtained by successive moving into the interior of the domain from its boundary, 
using, for instance, first-order upwind finite differences~\cite{Tsitsiklis,sethian1999fast}.
Among other popular methods, we should mention, first of all, the fast sweeping method 
\cite{tsai2003fast,Zhao}, which uses a Gauss--Seidel-style update strategy 
to progress across the domain. Recently (see, for example, \cite{jeong2008fast}), 
a fast iterative method for eikonal equations is actively developed using triangular \cite{fu2011fast} 
and tetrahedral \cite{fu2013fast} grids. Other modern variants of the fast marching method, 
which are adapted, in particular, to modern computing systems of parallel architecture, 
have been studied and compared, e.g., in~\cite{gomez2015fast}.

The second class of algorithms is associated with a transition from (\ref{1}), (\ref{2}) 
to a linear or nonlinear BVP for an elliptic equation  \cite{belyaev2015variational}. 
Instead of equation (\ref{1}), we (see~\cite{li2005level}) minimize the functional
\[
 J(u) = \frac{1}{2} \int_{\Omega} (|\nabla u| - 1)^2 d \bm x .
\]
It is possible to solve the BVP for the Euler--Lagrange equation for this functional, which has the form
\[
 \triangle u - \nabla \cdot \left (\frac{\nabla u}{|\nabla u|} \right ) = 0 ,
 \quad \bm x \in \Omega .   
\]
In~\cite{fares2002differential}, the computational algorithm is based on
solving the nonlinear boundary value problem for $v = 1/u$.
The solution of the problem (\ref{1}), (\ref{2}) can be related to the solution of 
the homogeneous Dirichlet problem for $p$-Laplacian:
\[
 \nabla \cdot (|\nabla u_p|^{p-2} \nabla u_p ) = - 1 ,
 \quad \bm x \in \Omega .   
\]  
In this case (see, e.g., \cite{bhattacharya1989limits,kawohl1990family}), we have
\[
 u_p(\bm x) \rightarrow u(\bm x) \ \mathrm{as} \ p \rightarrow \infty ,
 \quad \bm x \in \Omega . 
\]  
Thus, to find the solution of the problem (\ref{1}), (\ref{2}), we need to solve the nonlinear BVPs.

In our study, we focus on solving auxiliary boundary value problems for linear equations. 
This approach (see~\cite{gurumoorthy2009schrodinger,sethi2012schrodinger})is based on 
a connection between the nonlinear Hamilton--Jacobi equation and the linear Schrodinger equation.
Let $v_\alpha (\bm x)$ be the solution of the boundary value problem
\begin{equation}\label{3}
 - \alpha^2 \triangle v_\alpha + v_\alpha = 0,
 \quad \bm x \in \Omega , 
\end{equation} 
\begin{equation}\label{4}
 v_\alpha (\bm x) = 1,
 \quad \bm x \in \partial \Omega .  
\end{equation} 
Then, for $u_\alpha (\bm x) = - \alpha \ln (v_\alpha(\bm x)$, we have
$u_\alpha (\bm x) \rightarrow u(\bm x)$ as $\alpha  \rightarrow 0$.
A similar approach, where the auxiliary functions $v_\alpha (\bm x)$
are associated with the solution of the unsteady heat equation, is considered in the paper \cite{crane2013geodesics}.

In the present paper, we consider the eikonal equation in an anisotropic medium that is a more general variant 
in comparison with (\ref{1}). Using the transformation $u(\bm x) = - \alpha \ln (v(\bm x))$, 
the corresponding BVP of type (\ref{3}), (\ref{4}) is formulated for the new unknown quantity.
In our case, $\alpha \rightarrow 0$ and so, we have a singularly perturbed BVP
for the diffusion--reaction equation~\cite{roos2008robust,miller2012fitted}.
The numerical solution is based on using standard Lagrangian finite elements
\cite{brenner2008mathematical,larson2013finite}.
The main attention is paid to the monotonicity of the approximate solution for the auxiliary problem.

The paper is organized as follows. A boundary value problem for the eikonal equation in an anisotropic medium 
is formulated in Section 2. Its approximate solution is based on a transition to a singularly perturbed 
diffusion--convection equation. In Section 3, an approximation in space is constructed using Lagrangian finite elements 
and the main features of the problem solution are discussed. Numerical experiments on the accuracy 
of the approximate solution are presented in Section 5 for model two-dimensional problems. 
The results of the work are summarized in Section 5.

\section{Transformation of BVP for the eikonal equation in an anisotropic medium}

In a bounded polygonal domain $\Omega \subset R^m$, $m=1,2,3$ with the Lipschitz continuous boundary $\partial\Omega$, 
we search the solution of the BVP for the eikonal equation
\begin{equation}\label{5}
 E u = 1,
 \quad \bm x \in \Omega . 
\end{equation} 
Define the operator $E$ as
\begin{equation}\label{6}
  E u = \sum_{i=1}^{m} a_i^2 \left ( \frac{\partial u}{\partial x_i} \right )^2  
\end{equation} 
with the coefficients $a_i({\bm x}) > 0$.
The equation (\ref{5}) is supplemented with the homogeneous Dirichlet boundary condition
\begin{equation}\label{7}
  u({\bm x}) = 0,
  \quad {\bm x} \in \partial \Omega .
\end{equation}
The basic problems of numerical solving the boundary value problem (\ref{1})--(\ref{3}) result from the nonlinearity 
of the equation (see the operator $E$).

Similarly to~\cite{gurumoorthy2009schrodinger,sethi2012schrodinger}, we introduce the transformation
\begin{equation}\label{8}
 v_\alpha (\bm x) = \exp \left (- \frac{u_\alpha (\bm x) }{\alpha} \right ) 
\end{equation} 
with a numerical parameter $\alpha > 0$.
This type of transformation is widely used in studying differential equations with
quadratic nonlinearity $Eu$ (see, e.g., \cite{bicadze}).
 
Define the elliptic second-order operator $L$ by the relation 
\begin{equation}\label{9}
 L u = \sum_{i=1}^{m} \frac{\partial }{\partial x_i} \left ( a_i^2  \frac{\partial u}{\partial x_i} \right ) .  
\end{equation} 
For (\ref{8}), we have
\[
\begin{split}
 a_i^2 \frac{\partial v_\alpha}{\partial x_i} = &- \frac{1}{\alpha} \exp \left (- \frac{u_\alpha(\bm x) }{\alpha} \right ) 
 a_i^2 \frac{\partial u_\alpha}{\partial x_i} , \\
 \frac{\partial }{\partial x_i} \left ( a_i^2 \frac{\partial v_\alpha}{\partial x_i} \right ) =
 & - \frac{1}{\alpha} \exp \left (- \frac{u_\alpha(\bm x) }{\alpha} \right ) 
 \frac{\partial }{\partial x_i} \left ( a_i^2 \frac{\partial u_\alpha}{\partial x_i} \right ) \\
 & + \frac{1}{\alpha^2} \exp \left (- \frac{u_\alpha(\bm x) }{\alpha} \right ) 
 a_i^2 \left (\frac{\partial u_\alpha}{\partial x_i} \right )^2 .
\end{split}
\] 
By virtue of this, we obtain
\[
 \alpha^2 L v_\alpha - v_\alpha =  \exp \left (- \frac{u_\alpha(\bm x) }{\alpha} \right )
 (E u_\alpha - 1 - \alpha L u_\alpha ) .
\] 
Let $u_\alpha(\bm x)$ satisfies the equation
\begin{equation}\label{10}
 \alpha L u_\alpha - E u_\alpha  = - 1,
 \quad \bm x \in \Omega , 
\end{equation} 
and the boundary conditions
\begin{equation}\label{11}
  u_\alpha ({\bm x}) = 0,
  \quad {\bm x} \in \partial \Omega .
\end{equation} 
Under these conditions, for $v_\alpha(\bm x)$, we have the equation
\begin{equation}\label{12}
 \alpha^2 L v_\alpha - v_\alpha = 0,
 \quad \bm x \in \Omega . 
\end{equation} 
In view of (\ref{8}), from (\ref{11}), we obtain the following boundary condition:
\begin{equation}\label{13}
  v_\alpha ({\bm x}) = 1,
  \quad {\bm x} \in \partial \Omega .
\end{equation}

The equation (\ref{10}) can be treated as a regularization of the Hamilton--Jacobi equation via
the method of vanishing viscosity \cite{bookEvans}.
The boundary value problem (\ref{10}), (\ref{11}) produces an approximate solution of the
problem (\ref{5}), (\ref{6}) for small values of $\alpha$:
\[
 u_\alpha (\bm x) \rightarrow u(\bm x) \ \mathrm{as} \ \alpha  \rightarrow 0 ,
 \quad \bm x \in \Omega . 
\]
In this case, $u_\alpha (\bm x)$ is defined according to (\ref{8}) from the
solution of the linear boundary value problem (\ref{12}), (\ref{13}).

\section{Numerical implementation} 

An approximate solution of the BVP (\ref{5})--(\ref{7}) is represented (see~(\ref{8})) in the form
\begin{equation}\label{14}
 u_\alpha (\bm x) = - \alpha \ln (v_\alpha (\bm x)) , 
\end{equation} 
at a sufficiently low value of $\alpha$. In this case, $v_\alpha (\bm x)$ 
is defined as the solution of the BVP (\ref{12}), (\ref{13}). 
In the present work, the numerical implementation of this approach is carried out on the basis of
standard finite-element approximations \cite{brenner2008mathematical,larson2013finite}.
The main features of the computational algorithm result from the fact that the BVP
of diffusion--reaction  (\ref{12}), (\ref{13}) at small $\alpha$ is singularly perturbed, i.e.,
we have a small parameter at higher derivatives \cite{holmes2012introduction,verhulst2005methods}. 

Let us consider a standard quasi-uniform triangulation of the domain $\Omega$
into triangles in the 2D case or tetrahedra for 3D case. 
Let
\[
 V_0 = \{ v \in H^1(\Omega) \ | \ v(\bm x) = 0, \ \bm x \in \partial \Omega \} ,
\] 
\[
 V_1 = \{ v \in H^1(\Omega) \ | \ v(\bm x) = 1, \ \bm x \in \partial \Omega \} .
\] 
Denote by $V^h_0 \subset  V_0$ and $V^h_1 \subset  V_1$ the linear finite-element spaces.

For the BVP (\ref{12}), (\ref{13}), we put into the correspondence the variational problem of 
finding the numerical solution $y \in V^h_1$ from the conditions
\begin{equation}\label{15}
 a(y, v) = 0, 
\quad \forall v \in V^h_0 .
\end{equation} 
By (\ref{9}), for the bilinear form, we have
\[
 a(y,v) = \int_{\Omega} \left ( 
 \sum_{i=1}^{m} \alpha^2 a_i^2  \frac{\partial y}{\partial x_i}  \frac{\partial v}{\partial x_i} 
 + y v\right ) d \bm x .
\] 

The differential problem (\ref{12}), (\ref{13}) satisfies the maximum principle.
In particular, this guarantees the positiveness of the solution. More precisely 
(see, e.g.,~\cite{protter2012maximum,gilbarg2015elliptic}), 
for points inside the domain $\Omega$, we have
\[
 0 < v_\alpha (\bm x) < 1,
 \quad \bm x \in \Omega . 
\]
This the most important property must be also fulfilled for the solution of the discrete problem (\ref{15}):
\begin{equation}\label{16}
 0 < y(\bm x) < 1,
 \quad \bm x \in \Omega .  
\end{equation}
If (\ref{16}) holds, we speak of monotone approximations for the solution of the diffusion--reaction problem.

Even for regular boundary value problems, where the parameter $\alpha$ in (\ref{12}) is not small,
monotone approximations can be constructed using linear finite elements
with restrictions on the computational grid (Delaunay-type mesh, see, for instance, 
\cite{letniowski1992three,huang2011discrete}).
Additional restrictions appear (see, e.g.,~\cite{ciarlet1973maximum,brandts2008discrete})
on the magnitude of the reaction coefficient.
With respect to our problem (\ref{12}), (\ref{13}), for the grid step size, we have 
$h \leq \mathcal{O}(\alpha)$.

Restrictions on the grid due to the reaction coefficient can be removed.
The standard approach is related to the correction of approximations for the reaction coefficient based on  
the lumping procedure (see, e.g.,~\cite{Thomee2006}).

The standard approach to the solution of singularly perturbed diffusion--reaction problems 
(see~\cite{roos2008robust,miller2012fitted}) is based on using computational grids with
refinements in the vicinity of boundaries. A refinement of the grid is directly related to the value of 
the small parameter $\alpha$.

Another possibility to monotonize the solution of the problem (\ref{12}), (\ref{13}) at small values of $\alpha$ 
is the following approach. As noted in the paper \cite{cai2014natural}, for singularly perturbed
problems for the diffusion--convection equation, the use of finite-element approximations of higher order 
not only increases the accuracy of the approximate solution, but improves the monotonicity property as well.
It is interesting to check whether there is the same effect in the numerical solution of singularly perturbed 
problems for the diffusion--reaction equations.

\section{Numerical experiments} 

\begin{figure}[ht] 
  \begin{center}
    \begin{tikzpicture}
	\filldraw [color=blue!15] (0,0) rectangle +(4,6);
	\filldraw [color=blue!15] (4,0) rectangle +(4,4);
       	\draw [-,line width=1] (0,0) -- (0,6) -- (4,6) -- (4,4) -- (8,4) -- (8,0) -- (0,0);
       	\draw [->] (0,0) -- (0,7);
       	\draw [->] (0,0) -- (9,0);
       	\draw [dashed] (0,4) -- (4,4);
       	\draw [dashed] (4,0) -- (4,4);
       	\draw [-,color=red] (0,0) -- (4,4);
        \draw  (-0.35,-0.35) node {$0$}; 
        \draw  (4,-0.35) node {$1$}; 
        \draw  (8,-0.35) node {$2$}; 
        \draw  (-0.35,4) node {$1$}; 
        \draw  (-0.35,6) node {$1.5$}; 
        \draw  (-0.5,6.75) node {$x_2$}; 
        \draw  (8.75,-0.5) node {$x_1$}; 
    \end{tikzpicture}
    \caption{Computational domain} 
   \label{f-1}
  \end{center}
\end{figure}
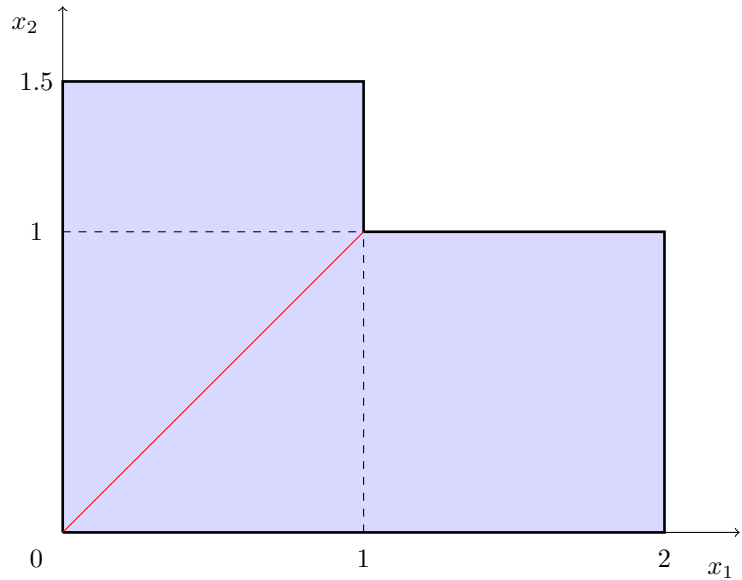

\begin{figure}[htp]
  \begin{center}
    \includegraphics[scale = 0.25] {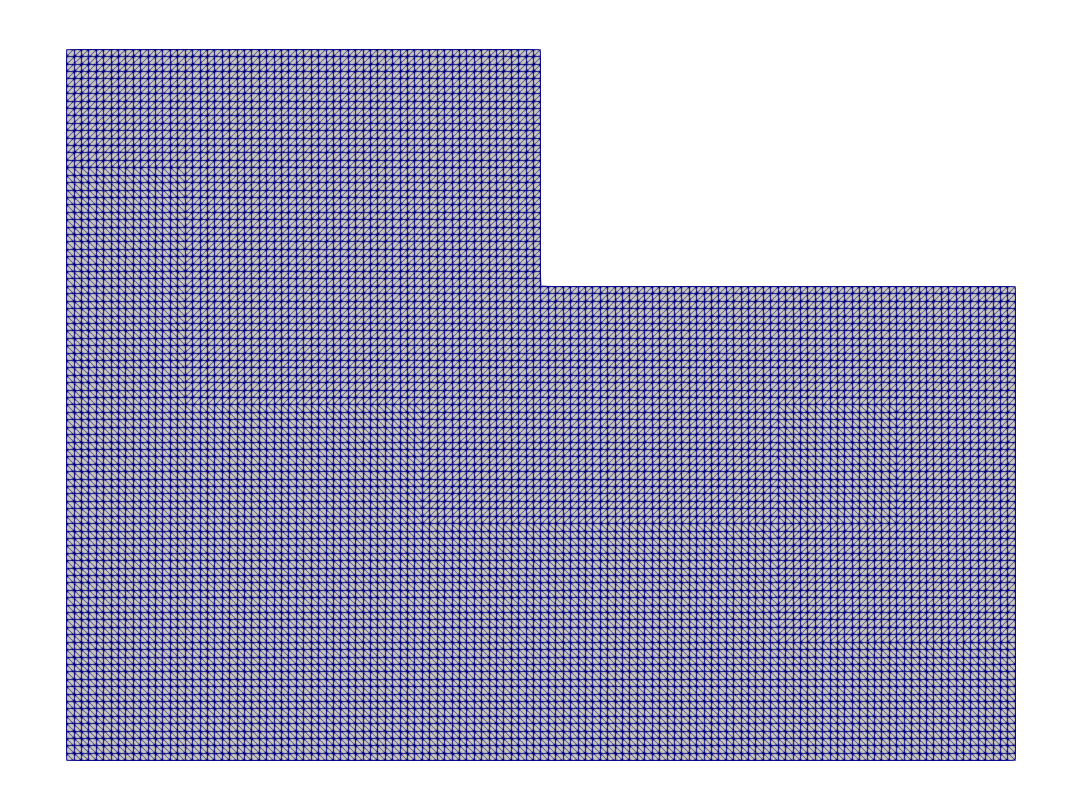}
	\caption{Basic (medium) computational grid}
	\label{f-2}
  \end{center}
\end{figure}

The 2D BVP (\ref{5})--(\ref{7}) in the L-shaped region depicted in Fig.~\ref{f-1} 
is considered as a model problem. We start with the simplest case, when $a_i = 1, \ i = 1, 2$.
The calculations have been performed on various grids. The basic (medium) computational grid, which
contains 10,465 nodes and 20,480 triangles, is shown in Fig.~\ref{f-2}.

\begin{figure}[htp]
  \begin{center}
    \includegraphics[scale = 0.3] {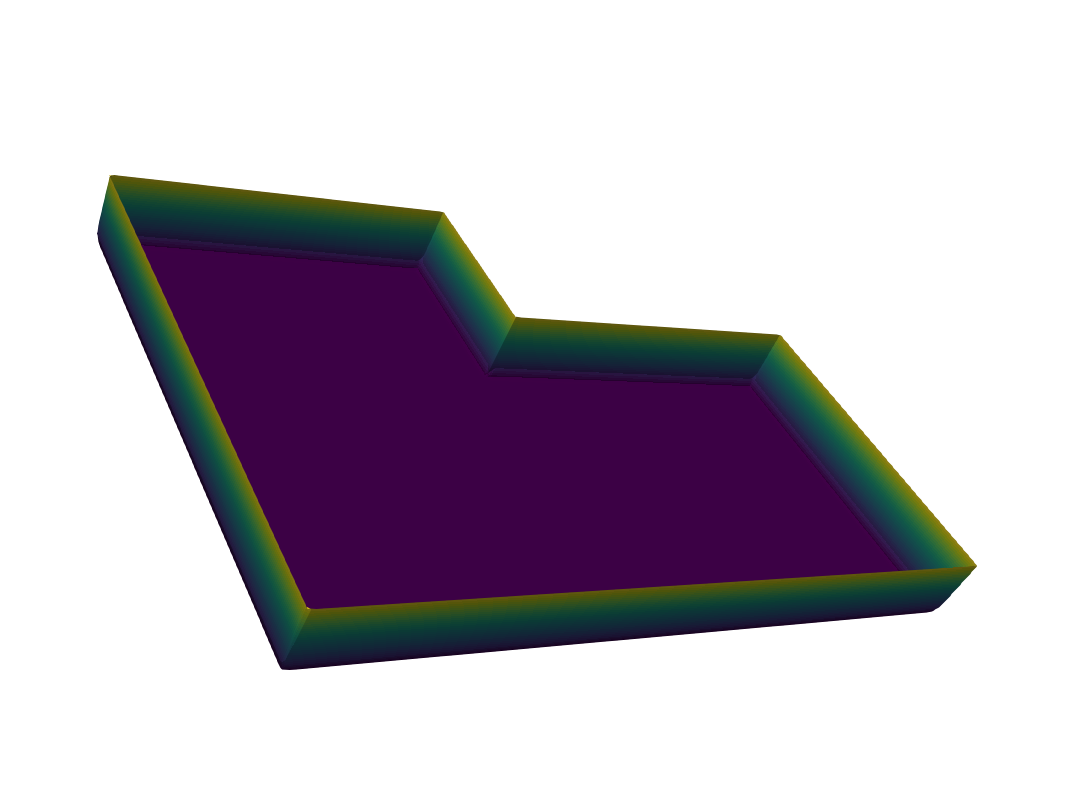}
	\caption{Solution $v_\alpha (\bm x)$ of the diffusion--reaction problem for $\alpha = 2^{-8}$}
	\label{f-3}
  \end{center}
\end{figure}

\begin{figure}[htp]
  \begin{center}
    \includegraphics[scale = 0.3] {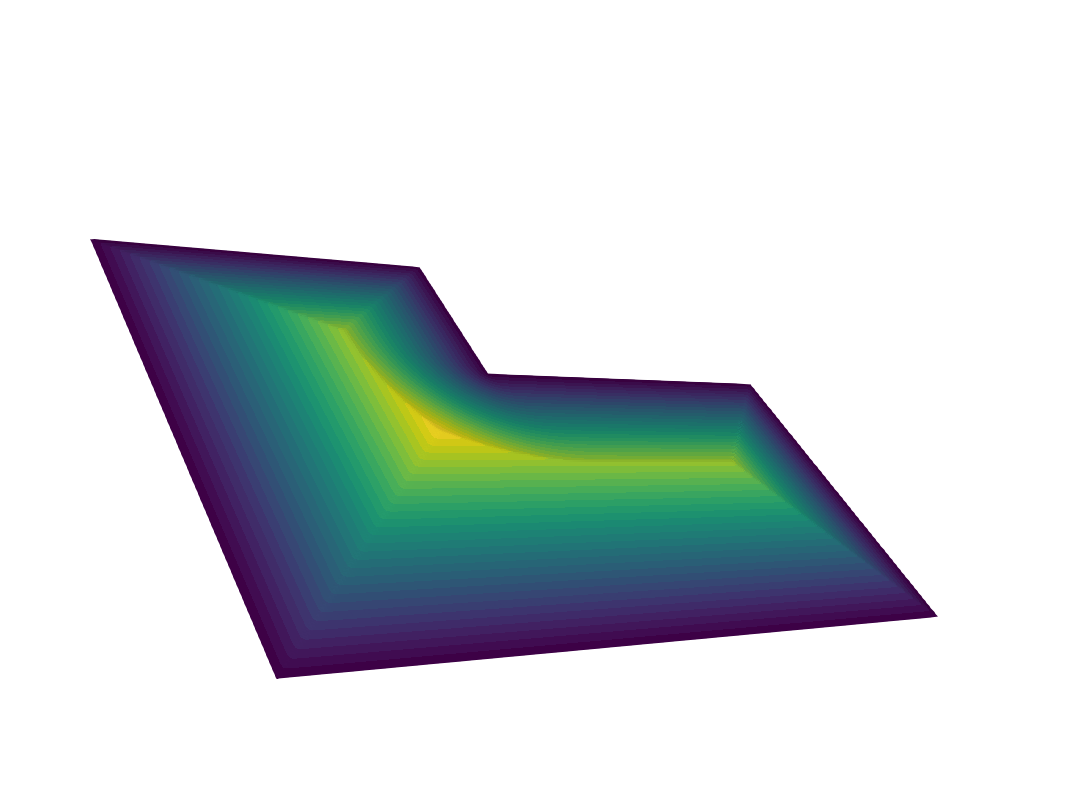}
	\caption{Solution $u_\alpha (\bm x)$ at $\alpha = 2^{-8}$}
	\label{f-4}
  \end{center}
\end{figure}

\begin{figure}[htp]
  \begin{center}
    \includegraphics[scale = 0.4] {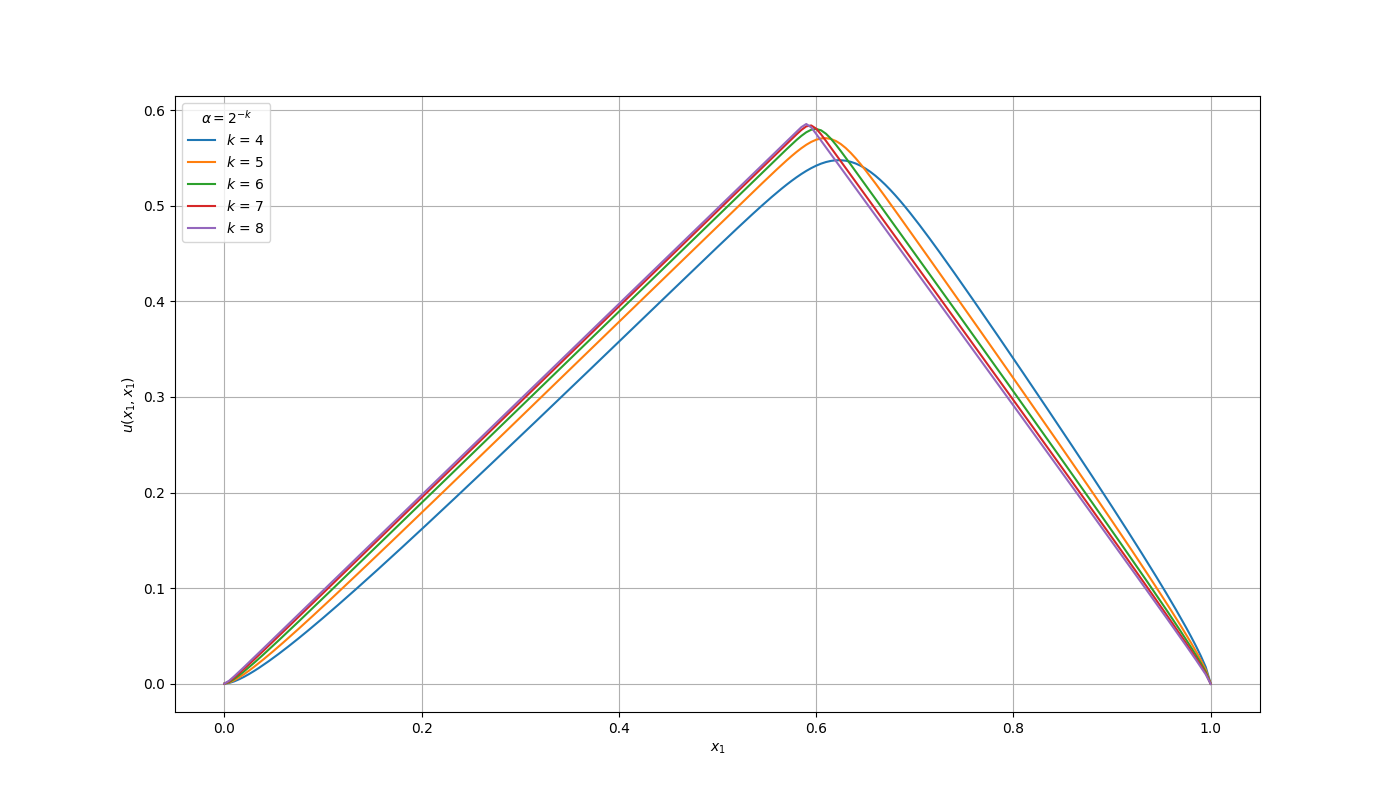}
	\caption{Solution $u_\alpha (\bm x)$ in the section $x_1 = x_2$ for various $\alpha$}
	\label{f-5}
  \end{center}
\end{figure}

In solving this problem, the key point is the dependence of the solution on the small parameter $\alpha$.
The numerical solution obtained on a very fine grid with $\alpha = 2^{-8}$ is treated as the exact one.
The solution $v_\alpha (\bm x)$ of the auxiliary problem (\ref{12}), (\ref{13}) under these conditions 
is presented in Fig.~\ref{f-3}, and the corresponding function $u_\alpha (\bm x)$, determined according 
to (\ref{14}), is shown in Fig.~\ref{f-4}. The influence of the parameter $\alpha$ can be observed 
in Fig.~\ref{f-5}, where the solution in the cross section $x_1 = x_2$ is plotted (the red line in Fig.~\ref{f-5}). 
In our model problem, a good accuracy is achieved for $\alpha \approx 2^{-7}$.

\begin{figure}[htp]
  \begin{center}
    \includegraphics[scale = 0.4] {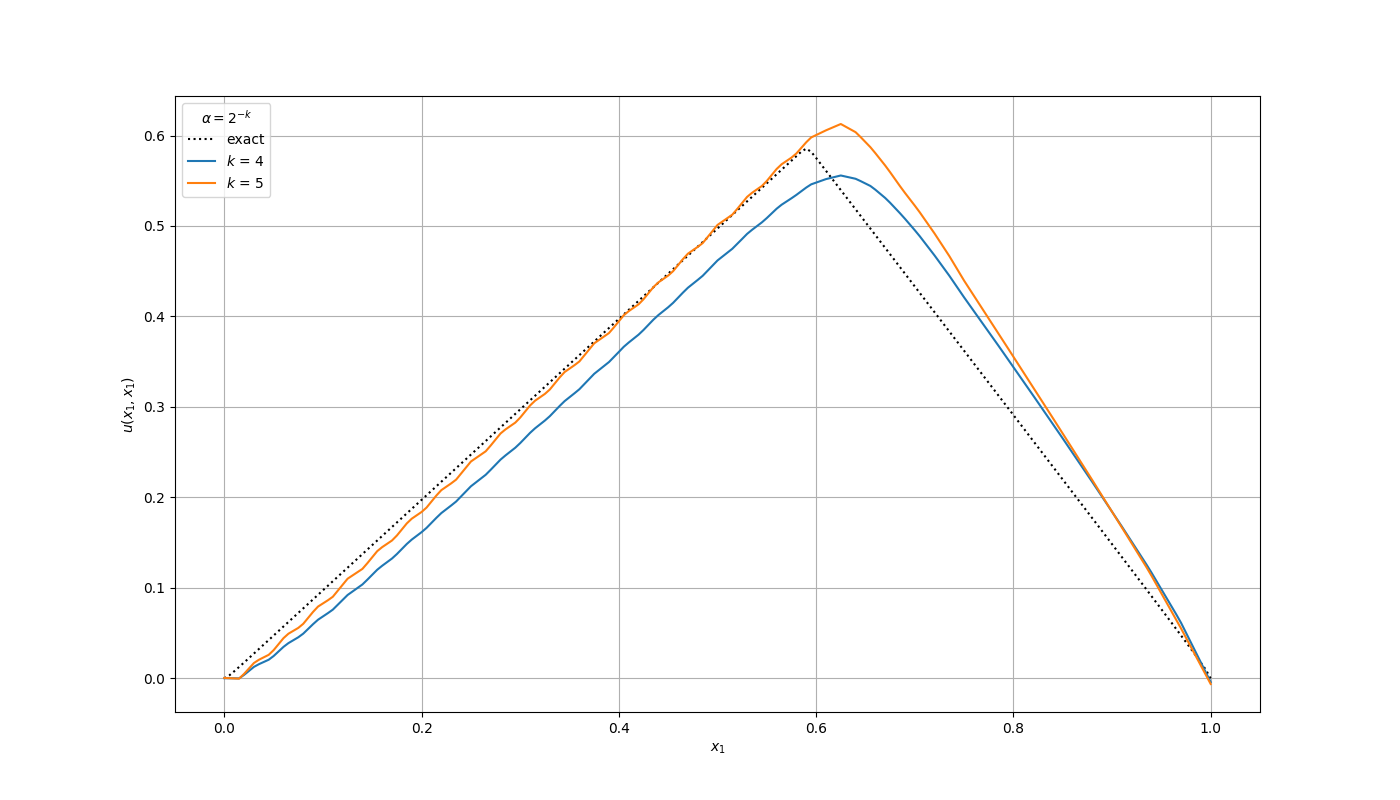}
	\caption{Solution $u_\alpha (\bm x)$ in the section $x_1 = x_2$ for various $\alpha$ --- the coarse grid}
	\label{f-6}
  \end{center}
\end{figure}

\begin{figure}[htp]
  \begin{center}
    \includegraphics[scale = 0.4] {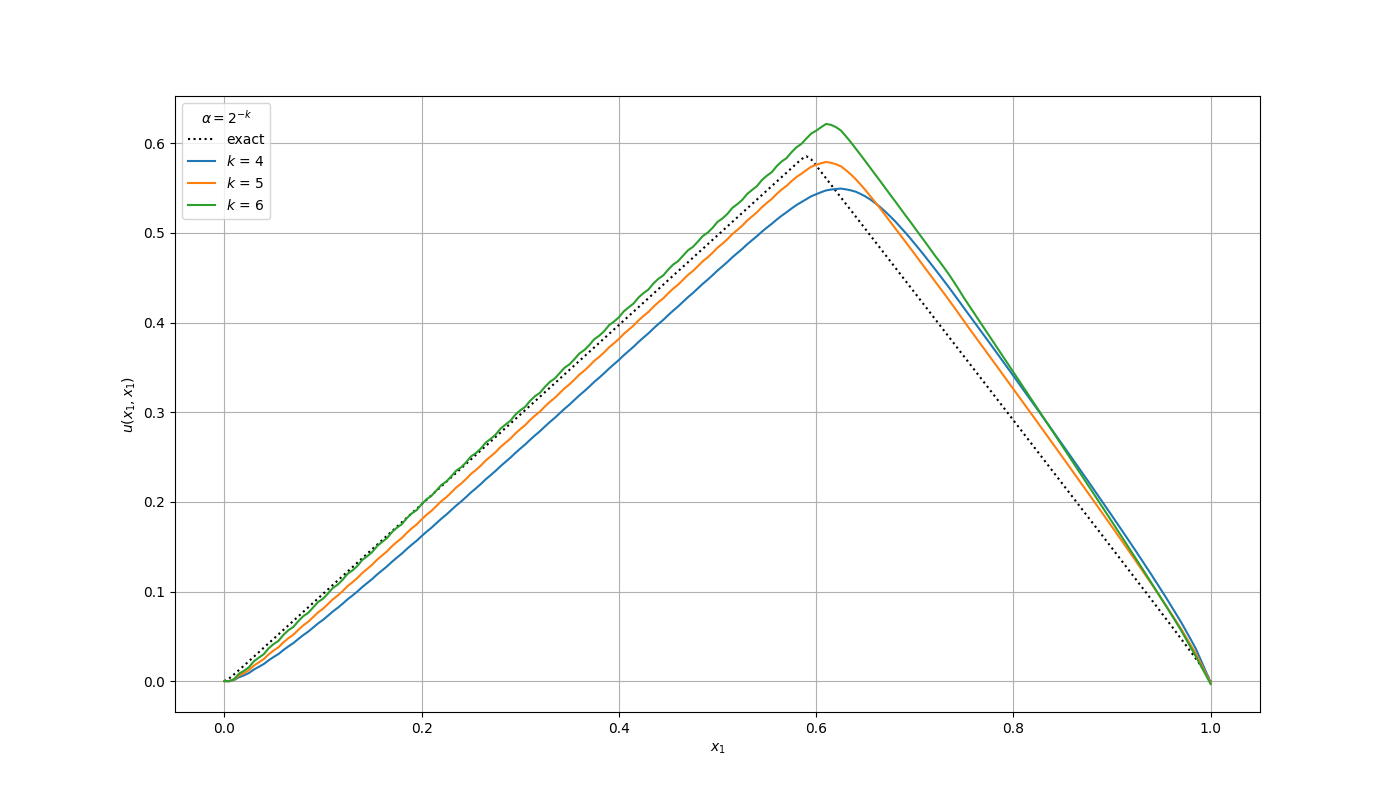}
	\caption{Solution $u_\alpha (\bm x)$ in the section $x_1 = x_2$ for various $\alpha$ --- the basic (medium) grid}
	\label{f-7}
  \end{center}
\end{figure}

\begin{figure}[htp]
  \begin{center}
    \includegraphics[scale = 0.4] {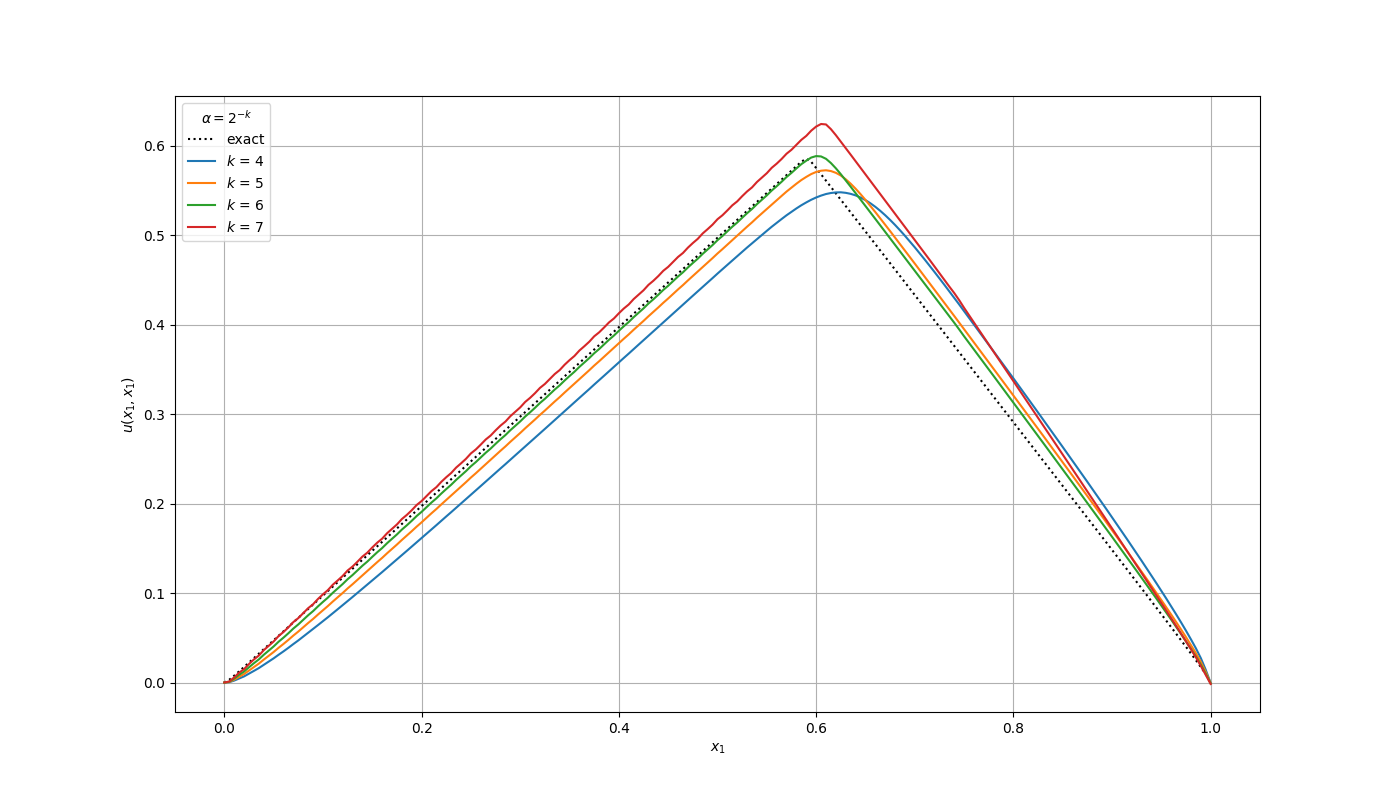}
	\caption{Solution $u_\alpha (\bm x)$ in the section $x_1 = x_2$ for various $\alpha$ --- the fine grid}
	\label{f-8}
  \end{center}
\end{figure}

The increase in accuracy can be achieved, first of all, by using finer grids. The solution for various $\alpha$ 
on the coarse grid (2,673 nodes and 5,120 triangles) is given in Fig.~\ref{f-6}. In this case, 
for $\alpha = 2^{-k}, \ k \geq  6$, the solution is non-monotone, i.e.,
at some nodes of the computational grid we have $y(\bm x) < 0$. Similar data for the basic grid are presented
in Fig.~\ref{f-7}. Here, the non-monotonicity appears at $\alpha = 2^{-k}, \ k \geq  7$.
Figure~\ref{f-8} demonstrates the numerical results obtained on the fine grid (41,409 nodes and 81,920 triangles). 
The non-monotonicity of the approximate solution occurs at $\alpha = 2^{-k}, \ k \geq  8$.

In the practical use of the approach (\ref{12})--(\ref{14}), it seems reasonable to follow the next strategy. 
We solve a number of auxiliary problems (\ref{12}), (\ref{13}) with a step-by-step decrease of 
the parameter $\alpha$ as long as the maximum principle holds. The solution obtained with the smallest $\alpha$ 
is taken as the approximate solution of the problem (\ref{5})--(\ref{7}). 

\begin{figure}[htp]
  \begin{center}
    \includegraphics[scale = 0.4] {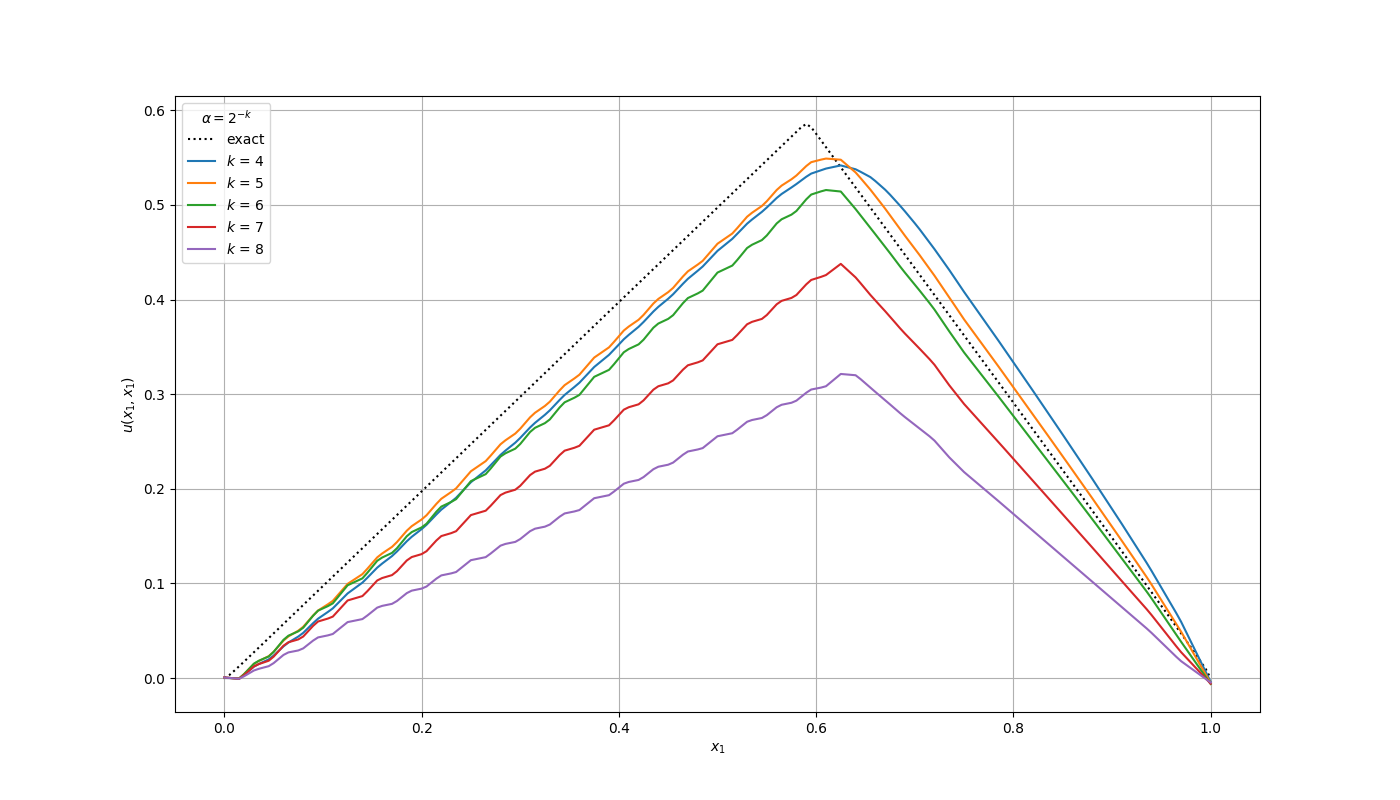}
	\caption{Reaction coefficient lumping for various $\alpha$ --- the coarse grid}
	\label{f-9}
  \end{center}
\end{figure}

\begin{figure}[htp]
  \begin{center}
    \includegraphics[scale = 0.4] {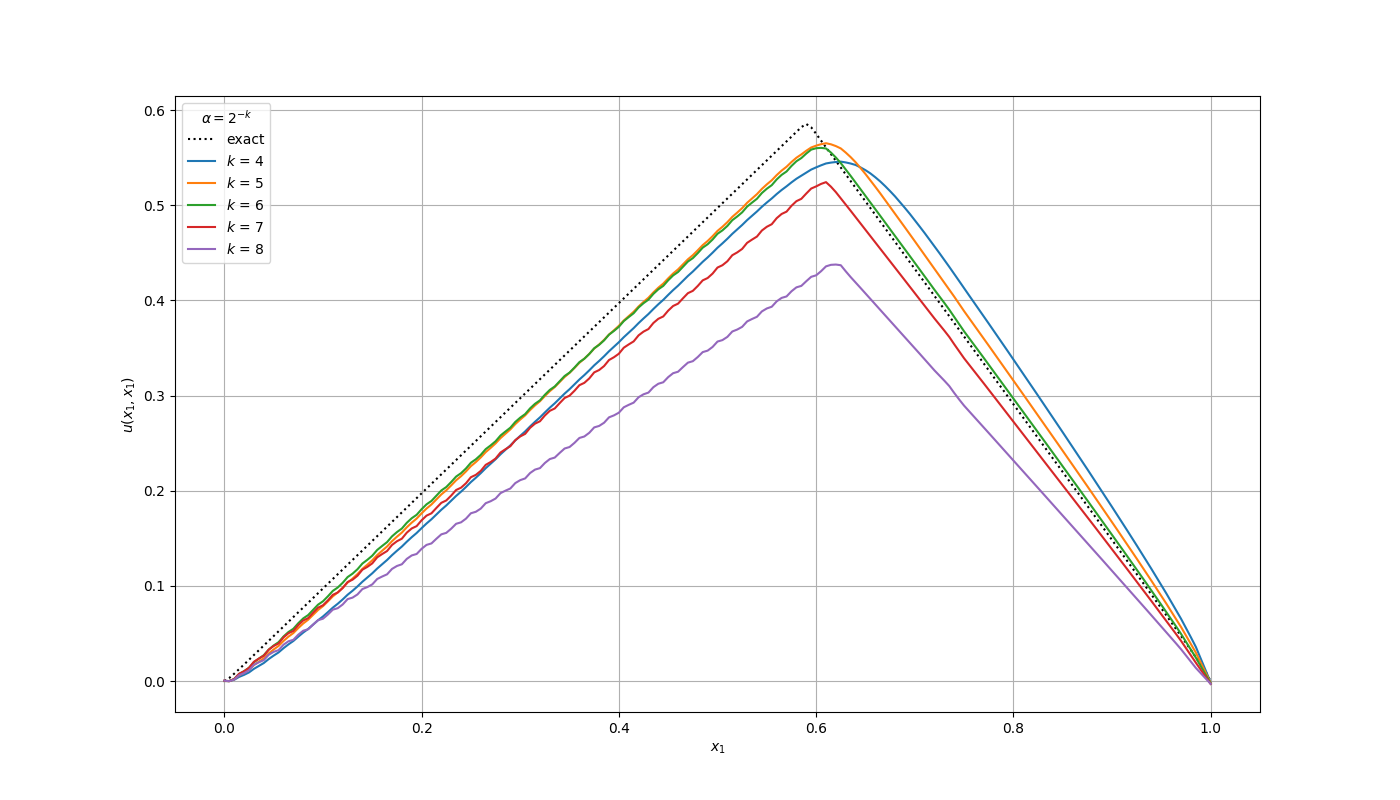}
	\caption{Reaction coefficient lumping for various $\alpha$ --- the basic (medium) grid}
	\label{f-10}
  \end{center}
\end{figure}

\begin{figure}[htp]
  \begin{center}
    \includegraphics[scale = 0.4] {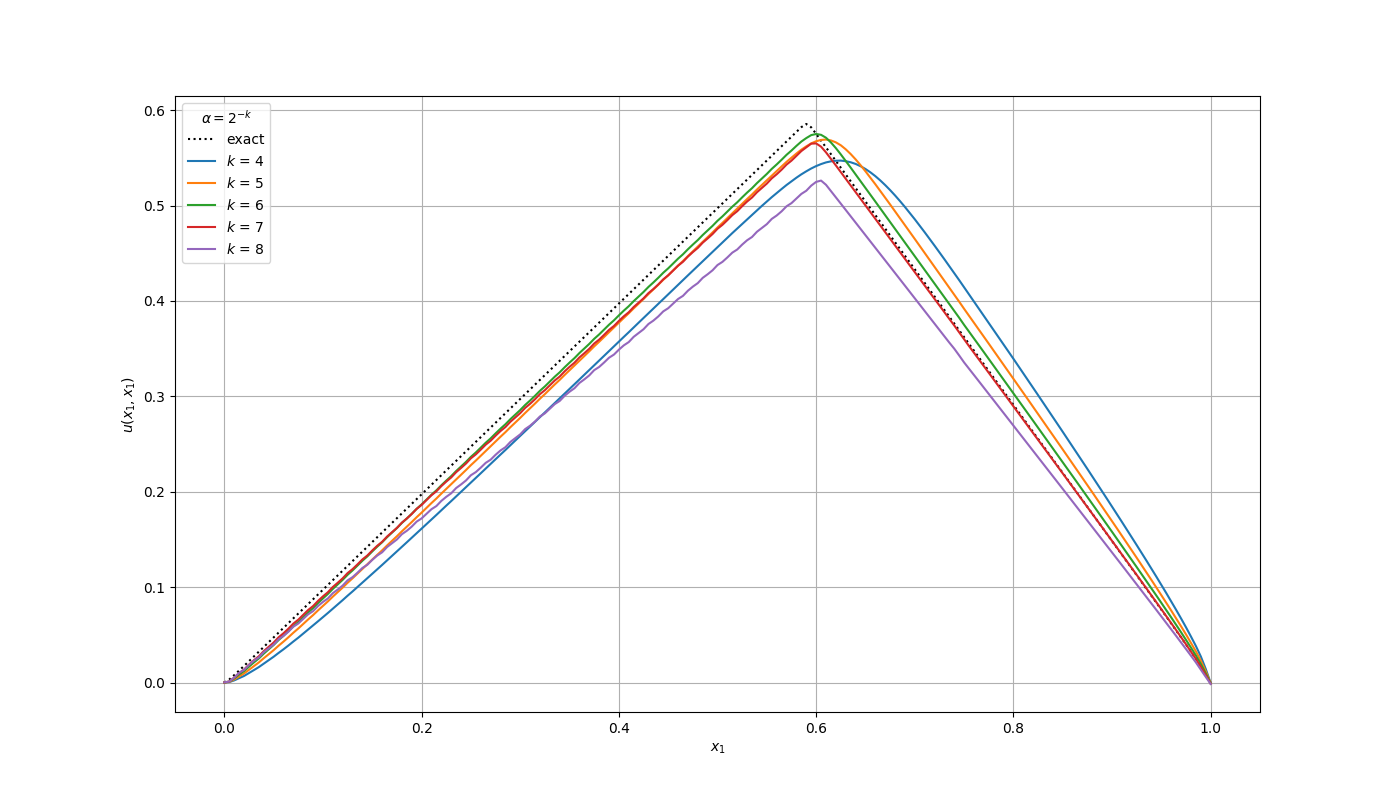}
	\caption{Reaction coefficient lumping for various $\alpha$ --- the fine grid}
	\label{f-11}
  \end{center}
\end{figure}

Our computational grids consist of rectangular isosceles triangles.
Because of this, the non-monotonicity is due to the reaction coefficient only.
To monotonize discrete solutions, it is sufficient to apply the standard procedure
of the reaction coefficient lumping \cite{Thomee2006}.
The effect of diagonalization of the reactive term in the finite-element approximation
in predictions on different computational grids can be observed in Figures~\ref{f-9}--\ref{f-11}.
In this case, the maximum principle holds for all $\alpha$.

The accuracy of the approximate solution decreases from some value of $\alpha$ as the parameter $\alpha$ decreases.
Moreover, the value of this optimal value is close to the value that we had without the lumping procedure. 
Therefore, we can use the diagonalization procedure for selecting the parameter $\alpha$ using the monotonicity condition 
for the discrete solution of the standard finite-element approximation.
In our case (see Figures~\ref{f-6}--\ref{f-8}), we select $\alpha=2^{-5}$ for the coarse grid, 
$\alpha=2^{-6}$ --- for the basic grid and $\alpha=2^{-7}$ --- for the fine grid.

\clearpage

\begin{figure}[htp]
  \begin{center}
    \includegraphics[scale = 0.4] {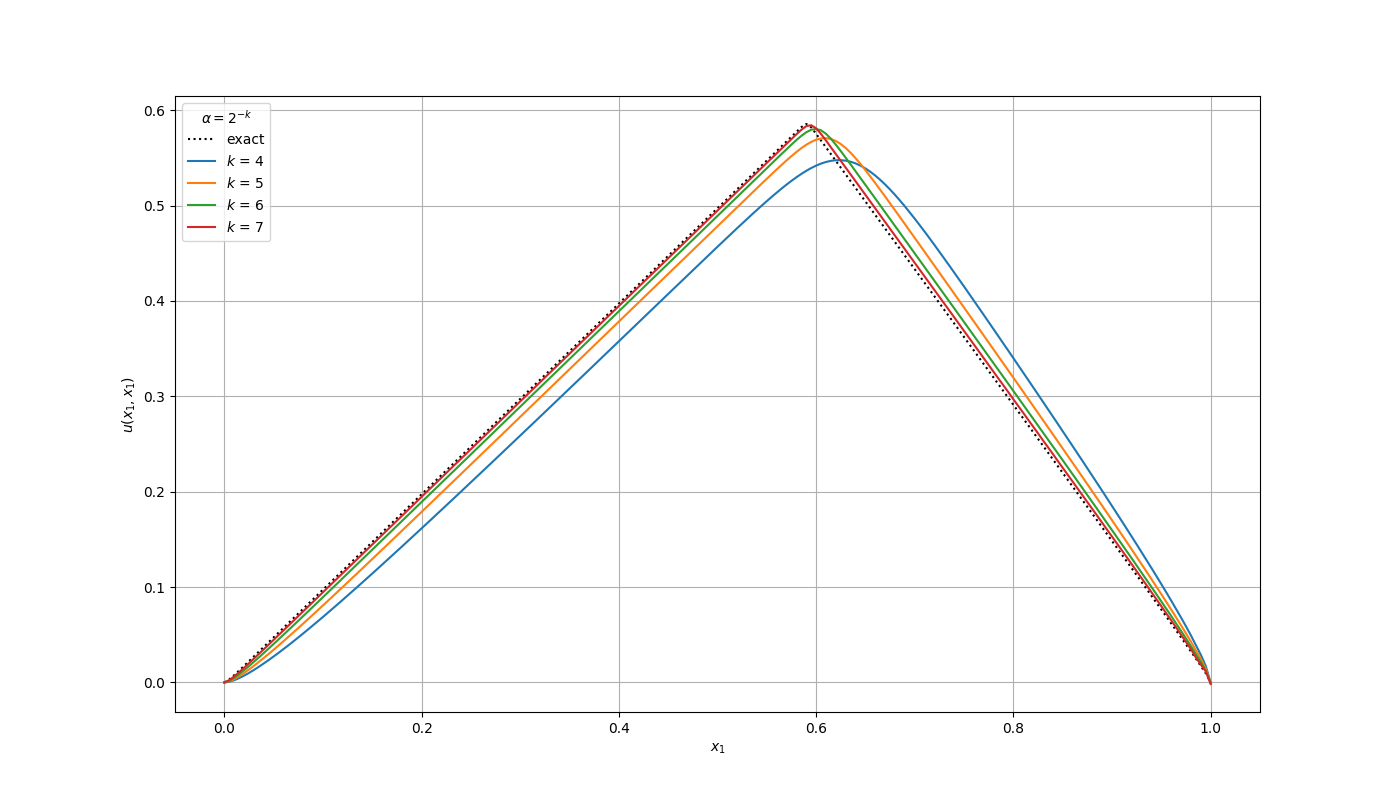}
	\caption{Solution $u_\alpha (\bm x)$ in the section  $x_1 = x_2$ for various $\alpha$ --- approximation $P_3$}
	\label{f-12}
  \end{center}
\end{figure}

\begin{figure}[htp]
  \begin{center}
    \includegraphics[scale = 0.4] {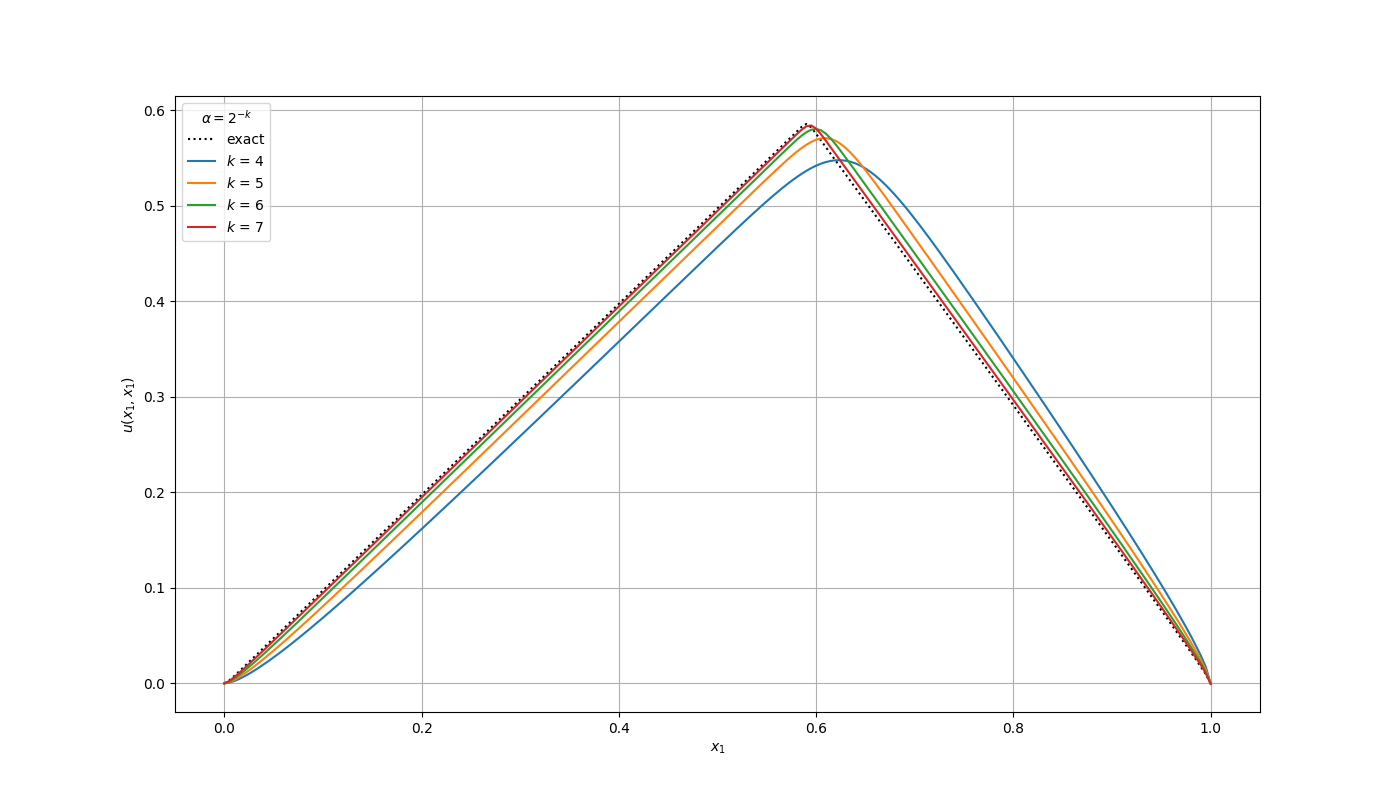}
	\caption{Solution $u_\alpha (\bm x)$ in the section  $x_1 = x_2$ for various $\alpha$  --- approximation $P_5$}
	\label{f-13}
  \end{center}
\end{figure}

\begin{figure}[htp]
  \begin{center}
    \includegraphics[scale = 0.4] {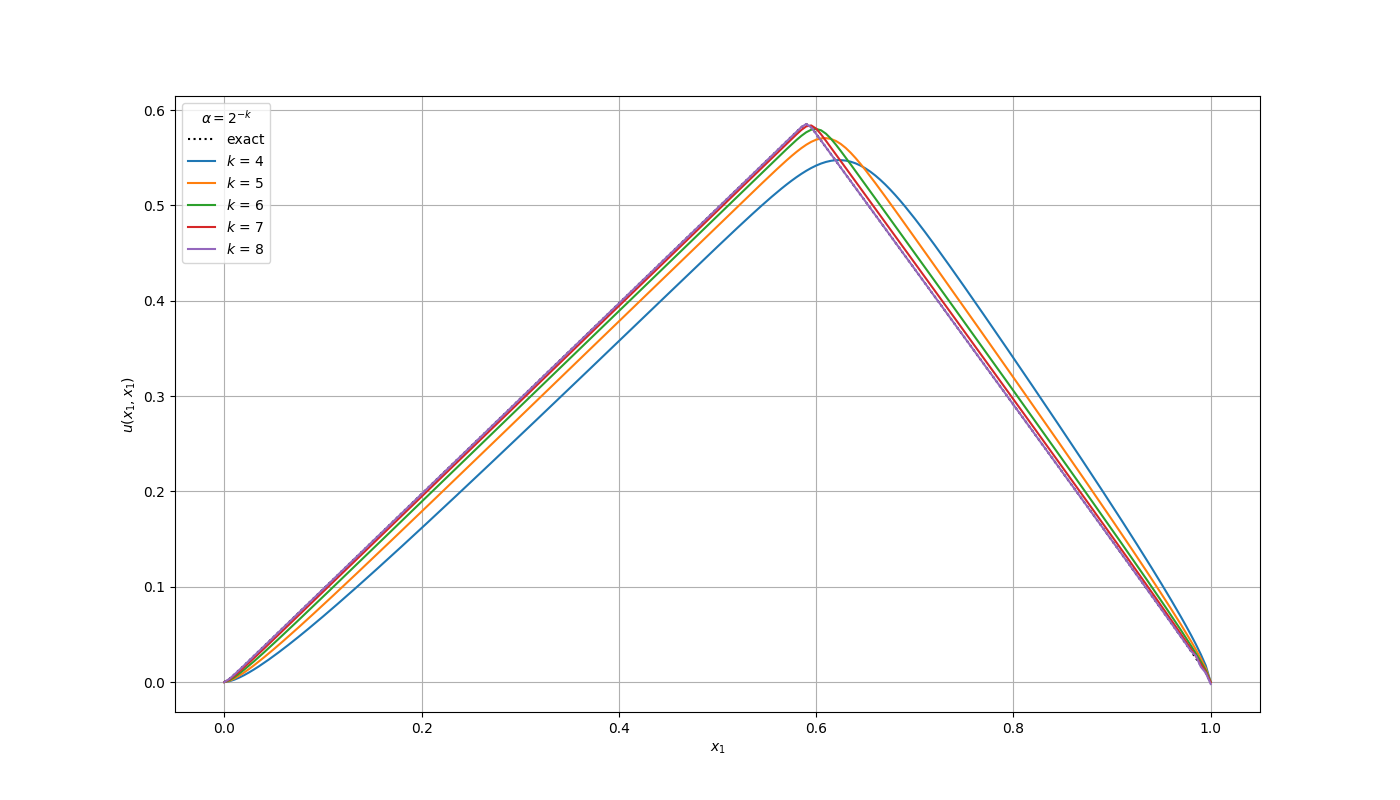}
	\caption{Solution $u_\alpha (\bm x)$ in the section  $x_1 = x_2$ for various $\alpha$  --- approximation $P_7$}
	\label{f-14}
  \end{center}
\end{figure}

Above, we have used linear finite elements. Below, we will present numerical results obtained on the basic grid
for Lagrangian finite elements of degree $m$, i.e., for approximations $P_m, \ m > 1$.
Figure~\ref{f-12} demonstrates the approximate solution obtained using finite elements of degree 3.
A comparison with the case $m=1$  (see Fig.~\ref{f-7}) indicates that the solution is more accurate and,
in addition, it is possible to carry out calculations with a smaller value of the parameter $\alpha$.
These effects become more pronounced when using finite elements of degree 5 (see Fig.~\ref{f-13}) and degree 7 
(see Fig.~\ref{f-14}).

Thus, the computational algorithm for solving the eikonal equation
(BVP (\ref{5})--(\ref{7})) can be based on the solution of the auxiliary problem (\ref{12}), (\ref{13}).
In doing so, we employ the minimum value of the parameter $\alpha$ that provides the monotone solution
on sufficiently fine computational grids using higher degree Lagrangian finite elements.

\clearpage

\begin{figure}[htp]
  \begin{center}
    \includegraphics[scale = 0.4] {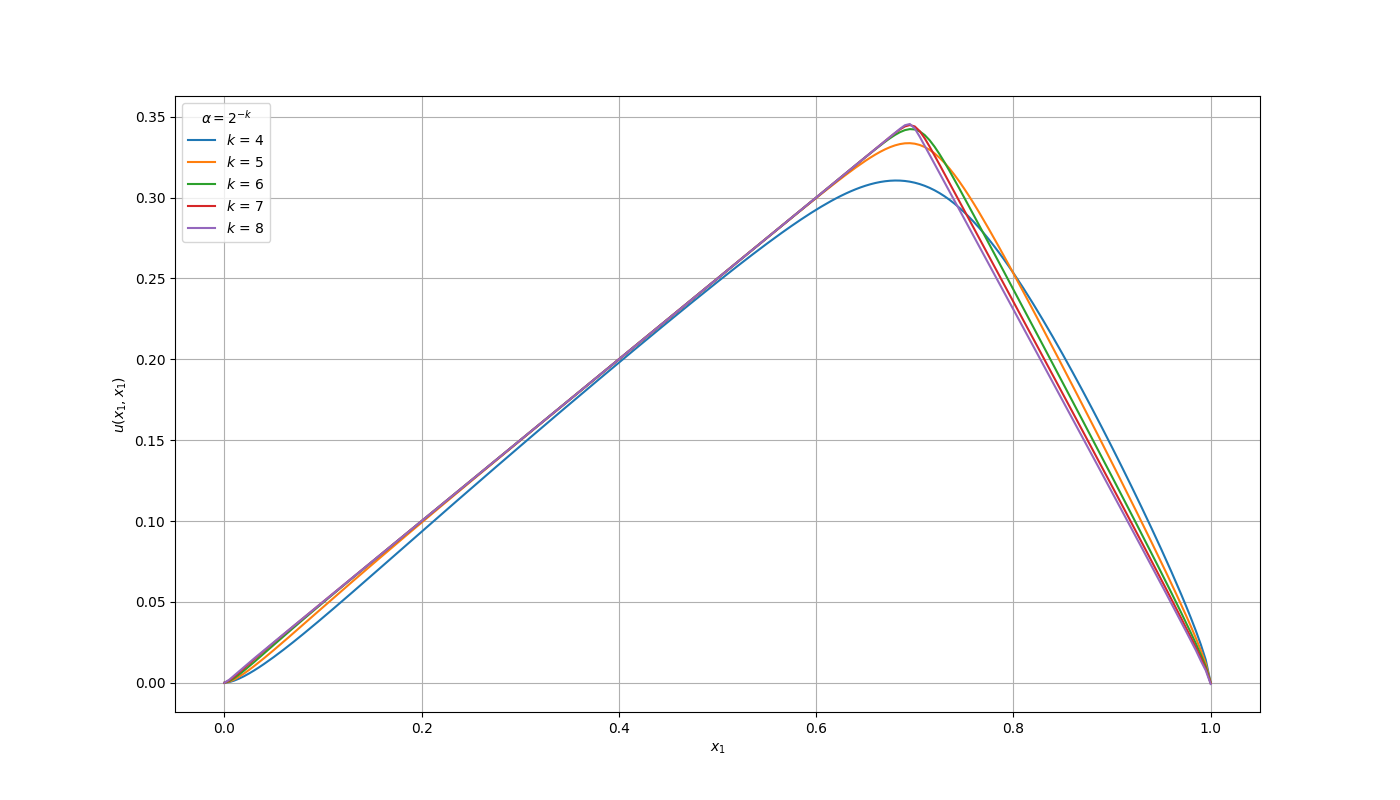}
	\caption{Solution $u_\alpha (\bm x)$ in the section $x_1 = x_2$ for various $\alpha$ --- $a_1^2 = 1, \ a_2^2 = 4$}
	\label{f-15}
  \end{center}
\end{figure}

\begin{figure}[htp]
  \begin{center}
    \includegraphics[scale = 0.3] {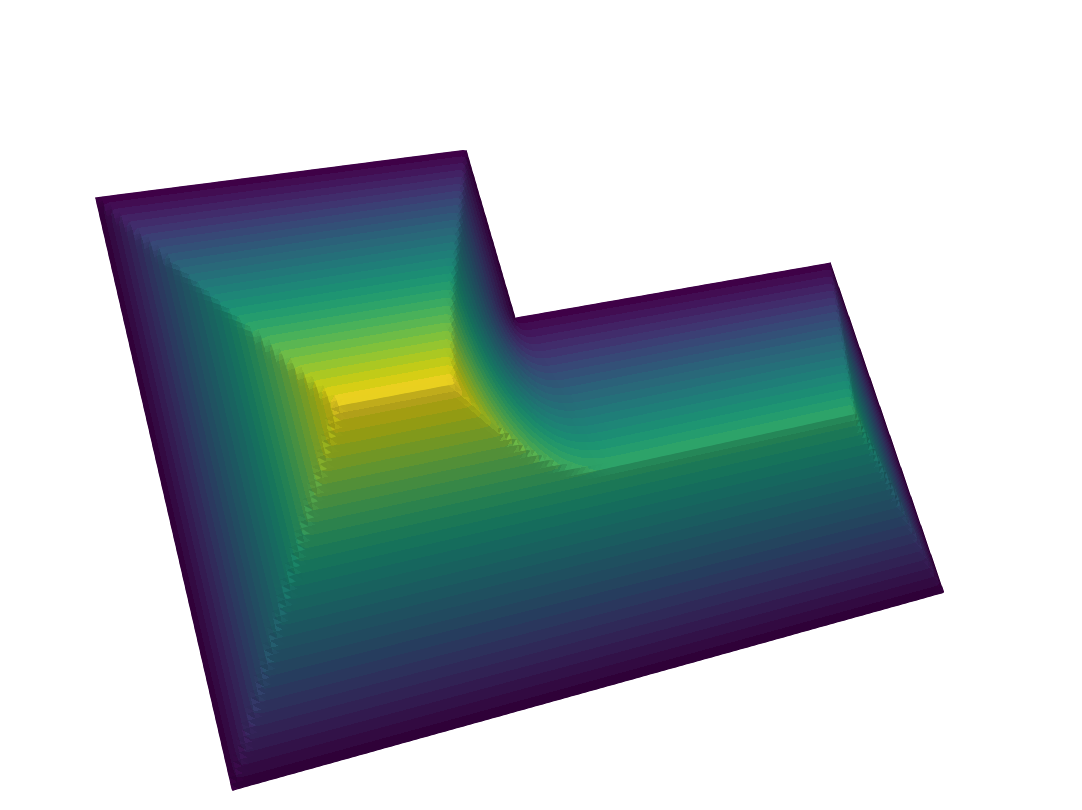}
	\caption{Approximate solution $u_\alpha (\bm x)$ for $\alpha= 2^{-8}$ at $a_1^2 = 1, \ a_2^2 = 4$}
	\label{f-16}
  \end{center}
\end{figure}

\begin{figure}[htp]
  \begin{center}
    \includegraphics[scale = 0.4] {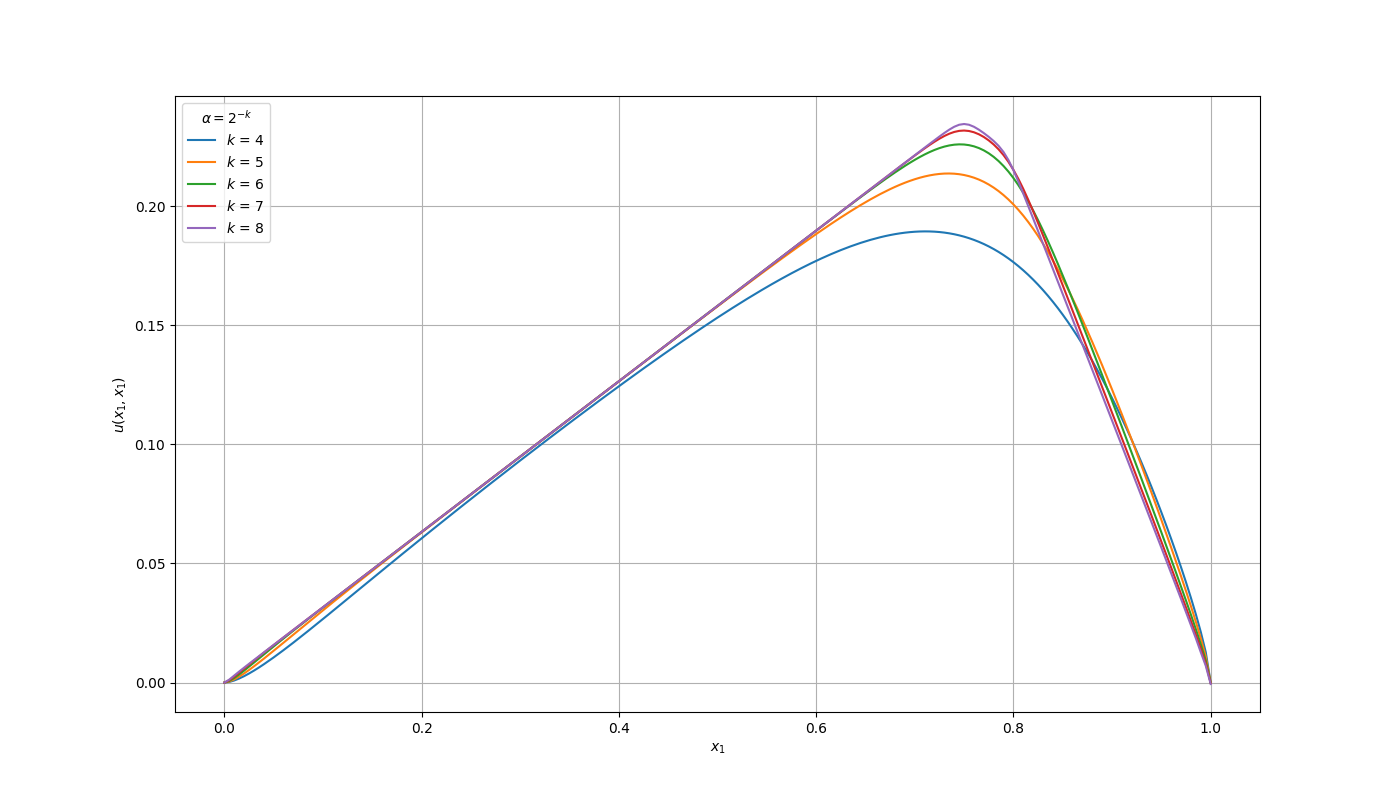}
	\caption{Solution $u_\alpha (\bm x)$ in the section $x_1 = x_2$ for various $\alpha$ --- $a_1^2 = 1, \ a_2^2 = 10$}
	\label{f-17}
  \end{center}
\end{figure}

\begin{figure}[htp]
  \begin{center}
    \includegraphics[scale = 0.3] {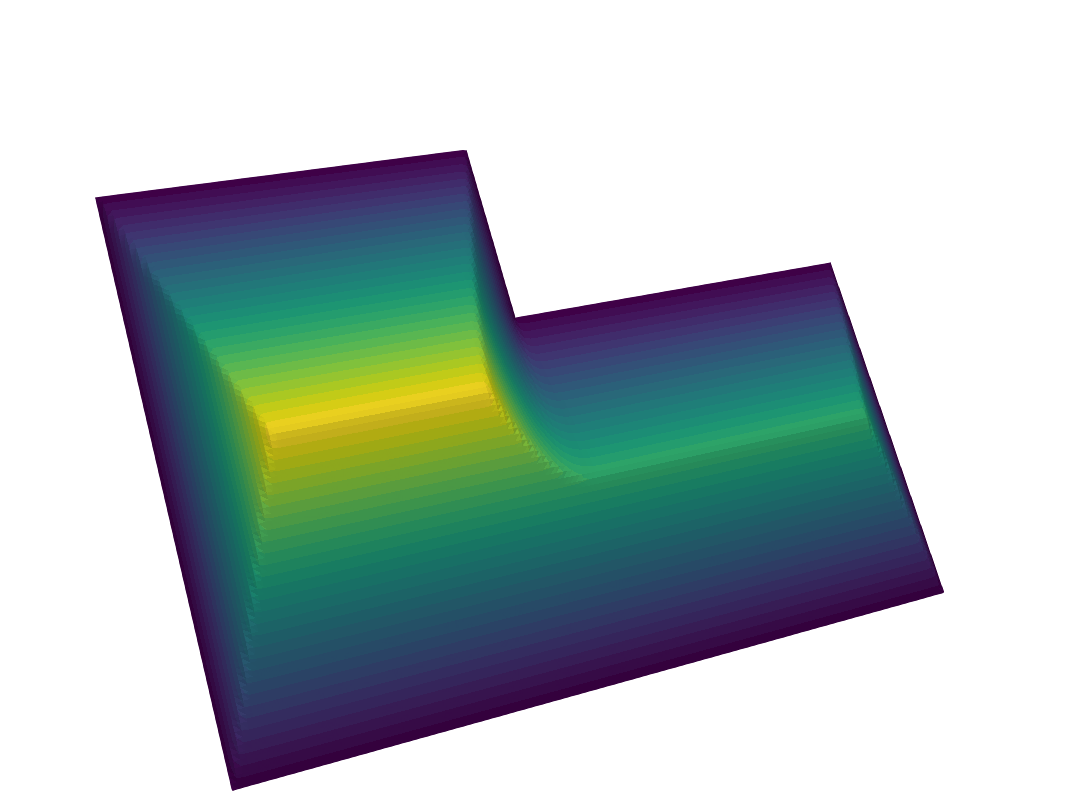}
	\caption{Approximate solution $u_\alpha (\bm x)$ for $\alpha= 2^{-8}$ at $a_1^2 = 1, \ a_2^2 = 10$}
	\label{f-18}
  \end{center}
\end{figure}

Special attention should be paid to the problem (\ref{5})--(\ref{7}) in the anisotropic case.
We have considered a variant with constant coefficients, where in (\ref{6}), we had
\[
 a_1^2 = 1,
 \quad a_2^2 = 4 . 
\]
The convergence of the approximate solution with decreasing $\alpha$ is given in Fig.~\ref{f-15}.
Calculations have been performed on the basic grid using finite-element approximation with $P_7$.
The numerical solution of the problem (\ref{5})--(\ref{7}) for $\alpha = 2^{-8}$ is shown in Fig.~\ref{f-16}. 
Similar data obtained for a more pronounced anisotropy:
\[
 a_1^2 = 1,
 \quad a_2^2 = 10, 
\]
are depicted in Fig.~\ref{f-17}, \ref{f-18}. 

\clearpage

\section{Conclusions} 

\begin{enumerate}
 \item A Dirichlet problem is considered for the multidimensional eikonal equation in a bounded domain with an anisotropic 
 medium.  The main peculiarities of such problems results from the fact that the eikonal equation is nonlinear.
 \item An approximate solution is constructed using a transformation of the original nonlinear boundary value problem 
 to a linear boundary value problem for the diffusion--reaction equation for an auxiliary function. 
 The transformed equation belongs to the class of singularly perturbed problems, i.e., there is a small parameter 
 at higher derivatives.
 \item Computational algorithms are constructed using standard finite-element approximations
 on triangular (2D problems) or tetrahedral (3D problems) grids.
 Monotonization of a discrete solution is achieved not only by using finer grids,
 but also via a correction of approximations for the reaction coefficient using the lumping procedure. 
 The use of finite elements of high degree is studied, too.
 \item Numerical experiments have been performed for 2D problems in order to demonstrate the robustness
 of the approach proposed in the work for solving boundary value problems for the eikonal equation in an anisotropic 
 medium. In particular, a good accuracy is observed when using sufficiently fine grids and Lagrangian finite elements 
 of higher degree.
\end{enumerate}

\section*{Acknowledgments}

Petr Vabishchevich gratefully acknowledges support from the
the Russian Federation Government (\#~14.Y26.31.0013).


\end{document}